\documentclass[11pt]{article}
\usepackage[a4paper,scale=0.8]{geometry}
\usepackage{amsmath}
\usepackage{amsthm}
\usepackage{amssymb}
\usepackage{bm}  
\usepackage{url}
\usepackage[numbers]{natbib}
\usepackage{enumerate}
\usepackage{graphicx}
\usepackage{blkarray}

\newcommand{\vct}{\bm}
\newcommand{\A}{{\mathrm{A}}}
\newcommand{\B}{{\mathrm{B}}}

\newcommand{\K}{{\mathrm{K}}}
\newcommand{\N}{{\mathrm{N}}}
\newcommand{\NP}{{\mathrm{NP}}}
\newcommand{\coNP}{{\mathrm{coNP}}}
\newcommand{\Tri}{\Delta}
\newcommand{\RR}{{\mathbb{R}}}
\newcommand{\trans}{{\mathrm{T}}}
\newcommand{\CUT}{{\mathrm{CUT}}}
\newcommand{\CutP}{\mathrm{CUT}^\square}
\newcommand{\symdiff}{\mathbin{\triangle}}
\newcommand*{\abs}[1]{\lvert#1\rvert}

\title{Generating facets for the cut polytope of a graph \\
  by triangular elimination}
\author{David Avis$^1$, Hiroshi Imai$^{2,3}$ and Tsuyoshi Ito$^{2,4}$%
  \thanks{Email: \texttt{tsuyoshi@is.s.u-tokyo.ac.jp}, Fax: +81-3-5800-6933.}}
\date{June 21, 2006}

\theoremstyle{plain}
\newtheorem{theorem}{Theorem}
\newtheorem{corollary}{Corollary}
\newtheorem{lemma}{Lemma}
\newtheorem{proposition}{Proposition}
\theoremstyle{definition}
\newtheorem{definition}{Definition}
\newtheorem{example}{Example}
\theoremstyle{remark}
\newtheorem{remark}{Remark}

\let\cite\citep

\BAnewcolumntype{s}{>{$\small$}c}

\begin{document}
\maketitle
\footnotetext[1]{School of Computer Science, McGill University,
  3480 University Street, Montreal, Quebec, Canada H3A 2A7.}
\footnotetext[2]{Department of Computer Science,
  Graduate School of Information Science and Technology,
  The University of Tokyo,
  7-3-1 Hongo, Bunkyo-ku, Tokyo, 113-0033 Japan.}
\footnotetext[3]{ERATO-SORST Quantum Computation and Information Project,
  Japan Science and Technology Agency,
  5-28-3 Hongo, Bunkyo-ku, Tokyo, 113-0033 Japan.}
\footnotetext[4]{Japan Society for the Promotion of Science.}

\begin{abstract}
  The cut polytope of a graph arises in many fields. Although much is known
  about facets of the cut polytope of the complete graph, very
  little is known for general graphs.
  The study of
  Bell inequalities in quantum information science requires knowledge
  of the facets of the cut polytope of the complete
  bipartite graph or, more generally, the complete $k$-partite graph.
  Lifting is a central tool to prove certain inequalities are facet
  inducing for the cut polytope.
  In this paper we introduce a lifting operation, named triangular
  elimination, applicable to the cut polytope of a wide range of
  graphs.
  Triangular elimination is a specific combination of zero-lifting and
  Fourier-Motzkin elimination using the triangle inequality.
  We prove sufficient conditions for the triangular elimination of
  facet inducing inequalities to be facet inducing.
  The proof is based on a variation of the lifting lemma adapted to
  general graphs.
  The result can be used to derive facet inducing inequalities of the
  cut polytope of various graphs from those of the complete graph.
  We also investigate the symmetry of facet inducing inequalities
  of the cut polytope of the complete bipartite graph derived by
  triangular elimination.
\end{abstract}

\section{Introduction}

\paragraph*{Cut polytope and related polytopes.}
The cut polytope arises in many fields
\cite{DezLau-JCAM94:applications1,DezLau-JCAM94:applications2,%
  DezLau:cut97},
and the structure of facets of the cut polytope has been intensively
studied.
For the complete graph with $n$ nodes,
a complete list of the facets of the cut polytope
$\CutP_n$ is known for
$n\le7$~\cite{Gri-EJC90},
as well as many classes of facet producing valid inequalities.
The hypermetric inequalities (see Chapter~28 of \cite{DezLau:cut97}) and
the clique-web inequalities~\cite{DezLau-MP92:facets2} (also Chapter~29 of
\cite{DezLau:cut97}), an extension of hypermetric
inequalities, are examples of such classes.
Very little is known about classes of facets for the cut polytope of
an arbitrary graph.
One such class are the cycle inequalities, which are projections of
the triangle inequalities.
They were shown to be facet producing by Barahona and
Majoub~\cite{BarMah-MP86}.
The structure of facets of the cut polytope is of both
theoretical and practical interest.
In the branch-and-cut approach to solve the MAX-CUT problem, facets of
the cut polytope are the most powerful cutting planes.
However, under the reasonable assumption that $\NP\ne\coNP$,
the complete list of the facets of the cut
polytope does not have a compact representation~\cite{Pit:prob89},
even for the complete graphs~\cite{AviDez-Net91,Pit-MP91}.
This implies that we cannot hope to enumerate all of its facets, but
rather should look for strong valid inequalities.

A lifting operation is a procedure which converts a given valid
inequality of the cut polytope of a small graph to a new valid
inequality of the cut polytope of a larger graph, and
is an established method for deriving new
facets systematically.
The most fundamental example of the lifting operations is zero-lifting
(for the complete graph~\cite{Dez-CRASParis73,DezLau-MP92:facets1} and
for general
graphs~\cite{Des-ORL90}).  Readers are referred to Section~26.5 and
Chapters~28--30 of \cite{DezLau:cut97} for more examples of classes of
valid inequalities and lifting operations.

The MAX-CUT problem is equivalent to unconstrained quadratic 0-1
programming~\cite{Ham-OR65}, and the associated boolean quadric polytope
is linearly isomorphic to the cut polytope.
This linear isomorphism is called the covariance mapping (see
Section~5.2 of \cite{DezLau:cut97}).
The boolean quadric polytope is also known as the correlation polytope,
especially in the physics literature.

\paragraph*{Relation of the cut polytope to quantum information processing.}
The polytopes described in the previous section have many applications in
quantum physics and quantum information theory~%
\cite{DezLau-JCAM94:applications2,DezLau:cut97}.
McRae and Davidson~\cite{McrDav-JMP72} showed the power of polytope
theory in quantum physics by proving that the possible solutions to
some problems arising in quantum physics form a convex polytope and
deriving inequalities for such solutions by convex hull algorithms.
One of the polytopes discussed there is identical to the boolean
quadric polytope.

In quantum information processing, the cut polytope and the boolean
quadric polytope arise in relation to Bell inequalities.
In this area, Bell inequalities, a generalization of Bell's original
inequality~\cite{Bel-Phys64},
are intensively studied~\cite{WerWol-QIC01,KruWer-0504166}
to better understand the nonlocality of quantum physics.
Bell inequalities deal with probabilities, and the search for
explicit formulae for Bell inequalities is related to Boole's
problem~\cite{Boo:laws1854}.
It is natural to consider Bell inequalities as inequalities valid for
certain convex polytopes~%
\cite{Fro-NC81,Pit-JMP86,Pit:prob89,Pit-MP91,Per:all99}
much in the same way as considering Boole's problem as a problem about
certain convex polytopes~\cite{DezLau-JCAM94:applications2}.
In particular, Bell inequalities involving joint probabilities of two
probabilistic events are exactly inequalities valid for the
boolean quadric polytope of a graph~\cite{Pit:prob89,Pit-MP91}.
To enumerate all the Bell inequalities for a given physical setting,
it is sufficient to enumerate the facets of the corresponding polytope by using
a convex hull algorithm.
Exhaustive enumeration of the Bell inequalities has been
performed~\cite{PitSvo-PRA01,ColGis-JPA04} in physical settings where
parameters such as the number of observables and the number of
possible outcomes of each observable are small enough.

Bell inequalities for two parties are inequalities valid for the
boolean quadric polytope of the complete bipartite graph $\K_{r,s}$,
and they correspond to inequalities valid for the cut polytope
$\CutP(\nabla\K_{r,s})$ via the covariance mapping.
$\nabla\K_{r,s}$ denotes the \emph{suspension graph} of $\K_{r,s}$,
that is, the graph obtained by adding a new node to $\K_{r,s}$ and
connecting it to all the existing nodes, and in other words, it is the
complete tripartite graph $\K_{1,r,s}$.
Enumeration of the facets of the cut polytope of the complete graph
uses symmetry and other structure specific to the cut polytope, and
they are often beyond the reach of general convex hull
packages.
Avis, Imai, Ito and Sasaki~\cite{AviImaItoSas-JPA05}
proposed an operation named \emph{triangular elimination}, which is a
combination of zero-lifting and Fourier-Motzkin elimination (see e.g.\
\cite{Zie:lectures98}) using the triangle inequality.
They proved that triangular elimination maps facet inducing inequalities of the
cut polytope of the complete graph to facet inducing inequalities of
the cut polytope of $\nabla\K_{r,s}$.

The cut polytope of $\nabla\K_{r,s}$ can be projected to the cut polytope
of $\K_{r,s}$, and this means that some Bell inequalities for the
correlation polytope of $\K_{r,s}$ correspond to inequalities valid
for $\CutP(\K_{r,s})$ via the covariance mapping.
Such Bell inequalities have good properties in relation to quantum
games~\cite{CleHoyTonWat-CCC04}.
They correspond to inequalities for correlation
functions~\cite{AviImaIto-0605148}, whose multi-party version is
discussed by Werner and Wolf~\cite{WerWol-PRA01} and \.{Z}ukowski and
Bruckner~\cite{ZukBru-PRL02}.

\paragraph*{Our results.}
In this paper, we generalize triangular elimination introduced in
\cite{AviImaItoSas-JPA05} to an operation which maps inequalities
valid for the cut polytope $\CutP(G)$ to those for $\CutP(G')$ for
graphs $G$ and $G'$ satisfying a certain condition.
From the viewpoint of combinatorial optimization, triangular
elimination is one of the lifting operations on inequalities valid for
the cut polytope.

Though the triangular elimination of an inequality is not uniquely
defined, all the choices are switching equivalent
(Proposition~\ref{prop:switching}) and therefore triangular elimination
can be seen as an operation which, given a switching equivalent class
of inequalities valid for $\CutP(G)$, uniquely produces a switching
equivalent class of inequalities valid for $\CutP(G')$.

We prove a sufficient condition (Theorem~\ref{theorem:trielim}) for
the triangular elimination of a facet inducing inequality to be facet
inducing.
The proof is similar to that of the zero-lifting theorem by Deza and
Laurent~\cite{DezLau-MP92:facets1,DezLau:cut97}, where the lifting
lemma used in the course of the proof is replaced with a version adapted
to general graphs.

For certain graphs $G$ and $G'$ which do not satisfy the conditions
of Theorem~\ref{theorem:trielim},
we can sometimes perform repeated triangular eliminations
on a sequence of graphs starting from $G$ and ending with $G'$.
Using this idea, we prove another sufficient condition in case where
$G=\K_n$.
This sufficient condition
extends Theorem~2.1 in \cite{AviImaItoSas-JPA05}.
It provides a method to derive a large number of inequalities which
define facets of the cut polytope of the complete $k$-partite graph.
These are relevant to $k$-party games, in light of the connection
between Bell inequalities and quantum games~\cite{CleHoyTonWat-CCC04}.

We also prove a necessary and sufficient condition for the triangular
eliminations of two facet inducing inequalities to be equivalent up to
permutation and switching in the case $G=\K_n$ and $G'=\K_{r,s}$.

\paragraph*{Organization of the paper.}
The rest of this paper is organized as follows.
Section~\ref{sect:prelim} reviews basic notions about the cut
polytope.
Section~\ref{sect:trielim} gives the definition of triangular
elimination for general graphs and proves its basic properties and the
main theorem stating a sufficient condition for the triangular
elimination of a facet to be a facet.
In Section~\ref{sect:k-partite}, we prove additional properties of
triangular elimination from the complete graph.
Section~\ref{sect:conclusion} states open problems.

\section{Preliminaries} \label{sect:prelim}

We briefly review basic notions about the cut polytope used in later
sections.
Definitions, theorems and other results stated in this section are
from the comprehensive reference~\cite{DezLau:cut97} on this topic,
which readers are referred to for more information.
We assume that readers are familiar with basic notions in convex
polytope theory such as convex polytope, facet, projection and
Fourier-Motzkin elimination.
Readers are referred to a textbook~\cite{Zie:lectures98} for details.

Throughout this paper, we use the following notation on graphs.
We denote the edge between two nodes $u$ and $v$ by $uv$.
For a graph $G=(V,E)$ and a node $v\in V$, we denote the neighbourhood
of $v$ by $\N_G(v)$.

\subsection{Cut polytope and cone}

The cut polytope (resp.\ cut cone) of a graph $G=(V,E)$ is the convex
hull (resp.\ conic hull) of the cut vectors of $G$.
A formal definition is as follows.

\begin{definition}[Cut polyhedra]
  The \emph{cut polytope} of a graph $G=(V,E)$, denoted $\CutP(G)$, is
  the convex hull of the cut vectors $\vct{\delta}_G(S)$ of $G$
  defined by all the subsets $S\subseteq V$ in the
  $\abs{E}$-dimensional vector space $\RR^E$.
  The cut vector $\vct{\delta}_G(S)$ of $G$ defined by $S\subseteq V$
  is a vector in $\RR^E$ whose $uv$-coordinate is defined as follows:
  \[
    \delta_{uv}(S)=\begin{cases}
      1 & \text{if $\abs{S\cap\{u,v\}}=1$,} \\
      0 & \text{otherwise,}
    \end{cases} \text{ for $uv\in E$.}
  \]
  The \emph{cut cone} of $G$, denoted $\CUT(G)$, is the conic hull of
  the cut vectors $\vct{\delta}_G(S)\in\RR^E$ of $G$ for all the
  subsets $S\subseteq V$.
  If $G$ is the complete graph $\K_n$, we denote $\CutP(\K_n)$ and
  $\CUT(\K_n)$ also as $\CutP_n$ and $\CUT_n$, respectively.
\end{definition}

For a subset $F$ of a set $E$, the \emph{incidence vector}
of $F$ (in $E$)%
\footnote{The set $E$ is sometimes not specified explicitly
  when $E$ is clear from the context or the choice of $E$ does not
  make any difference.}
is the vector $\vct{x}\in\{0,1\}^E$ defined by
$x_e=1$ for $e\in F$ and $x_e=0$ for $e\in E\setminus F$.
Using this term, the definition of the cut vector can also be stated
as follows: $\vct{\delta}_G(S)$ is the incidence vector of the cut set
$\{uv\in E\mid\abs{S\cap\{u,v\}}=1\}$ in $E$.

The cut polytope and cone are full-dimensional in
$\RR^E$~\cite{BarGroMah-MOR85}.
The following inequalities are the first class of facets of the cut
cone of an arbitrary graph.

\begin{theorem}[\cite{BarMah-MP86}] \label{theorem:cycle}
  \begin{enumerate}[(i)]
  \item
    For a graph $G=(V,E)$, a cycle $C\subseteq E$ in $G$ and an edge
    $uw\in C$, the \emph{cycle inequality}
    \begin{equation}
      x_{uw}-\sum_{e\in C\setminus\{uw\}}x_e\le0
      \label{eq:cycle}
    \end{equation}
    is valid for $\CUT(G)$.
  \item
    If $C$ is a chordless cycle in $G$, then (\ref{eq:cycle}) is facet
    inducing for $\CUT(G)$.
  \end{enumerate}
\end{theorem}

The following proposition follows immediately from the fact that the origin
is a vertex of $\CutP(G)$.

\begin{proposition} \label{prop:cone-homogeneous}
  Inequality $\vct{a}^\trans\vct{x}\le0$ is valid (resp.\ facet
  inducing) for $\CutP(G)$ if and only if it is valid (resp.\ facet
  inducing) for $\CUT(G)$.
\end{proposition}

\subsection{Operations on inequalities}

\subsubsection{Symmetric transformations}

Let $G=(V,E)$ be a graph.
The cut polytope $\CutP(G)$ admits two kinds of symmetric
transformations, which correspond to operations on valid inequalities
which preserve their properties.

\begin{definition}[Permutation]
  Let $\sigma$ a permutation on $V$ which is an automorphism of $G$.
  Then the \emph{$\sigma$-permutation} of an inequality
  $\vct{a}^\trans\vct{x}\le a_0$ is an inequality
  $(\vct{a}')^\trans\vct{x}\le a_0$ where $\vct{a}'\in\RR^E$ is
  defined by $a'_{ij}=a_{\sigma(i)\sigma(j)}$.
  Such an inequality is said to be \emph{permutation equivalent} to
  $\vct{a}^\trans\vct{x}\le a_0$.
\end{definition}

\begin{definition}[Switching]
  Let $S$ be a subset of $V$.
  Then the \emph{$S$-switching} of an inequality
  $\vct{a}^\trans\vct{x}\le a_0$ is an inequality
  $(\vct{a}')^\trans\vct{x}\le a_0-\vct{a}^\trans\vct{\delta}_G(S)$
  where $\vct{a}'\in\RR^E$ is defined by
  $a'_{ij}=(-1)^{\delta_{ij}(S)}a_{ij}$.
  Such an inequality is said to be \emph{switching equivalent} to
  $\vct{a}^\trans\vct{x}\le a_0$.
\end{definition}

Generalizing the cycle inequality in the form of (\ref{eq:cycle}),
the cycle inequality~\cite{BarMah-MP86} for the cut polytope
$\CutP(G)$ is defined as follows.
For a cycle $C\subseteq E$ in $G$ and a subset $F\subseteq C$ with
$\abs{F}$ odd,
\begin{equation}
  \sum_{e\in F}x_e-\sum_{e\in C\setminus F}x_e\le\abs{F}-1.
  \label{eq:cycle-sw}
\end{equation}
Inequality (\ref{eq:cycle-sw}) is
switching equivalent to (\ref{eq:cycle}), since
it is the $S$-switching of (\ref{eq:cycle}) where $S$ is a subset of
the nodes in $C$ such that the intersection of $C$ and the cut set
defined by $S$ is equal to $F\symdiff\{uw\}$.
Here $F\symdiff\{uw\}$ denotes the symmetric difference of the two
sets $F$ and $\{uw\}$.

We say $(\vct{a}')^\trans\vct{x}\le a'_0$ is
\emph{permutation-switching equivalent} to
$\vct{a}^\trans\vct{x}\le a_0$ if they can be transformed to each
other by using permutation and/or switching equivalence.

The following proposition is stated as Lemma~26.2.1 and
Corollary~26.3.7 in \cite{DezLau:cut97}.

\begin{proposition} \label{prop:symmetry}
  Let $\vct{a}^\trans\vct{x}\le a_0$ and
  $(\vct{a'})^\trans\vct{x}\le a'_0$ be permutation-switching
  equivalent inequalities.
  Then $\vct{a}^\trans\vct{x}\le a_0$ is valid (resp.\ facet inducing)
  for $\CutP(G)$ if and only if $(\vct{a}')^\trans\vct{x}\le a'_0$ is
  valid (resp.\ facet inducing) for $\CutP(G)$.
\end{proposition}

A \emph{root} of an inequality is a cut vector that satisfies it as
an equation.
The following well-known proposition, which follows from the definition
of switching, shows the essential equivalence of the cut cone and polytope.

\begin{proposition} \label{prop:cone-switching}
  Let $\vct{a}^\trans\vct{x}\le a_0$ be an inequality, valid for
  $\CutP(G)$, which has a root $\vct{\delta}_G(S)$.
  Then its $S$-switching $(\vct{a}')^\trans\vct{x}\le0$
  is valid for $\CUT(G)$.
\end{proposition}

From Theorem~\ref{theorem:cycle} and Proposition~\ref{prop:symmetry},
the following corollary follows immediately.

\begin{corollary}[\cite{BarMah-MP86}] \label{cor:cycle-sw}
  \begin{enumerate}
  \item
    For a graph $G=(V,E)$, a cycle $C\subseteq E$ in $G$ and a subset
    $F\subseteq C$ with $\abs{F}$ odd, the cycle
    inequality~(\ref{eq:cycle-sw}) is valid for $\CutP(G)$.
  \item
    If $C$ is a chordless cycle, then (\ref{eq:cycle-sw}) is facet
    inducing for $\CutP(G)$.
  \end{enumerate}
\end{corollary}

\subsubsection{Collapsing}

Let $uv$ be an edge of a graph $G=(V,E)$.
The intersection of the cut polytope $\CutP(G)$ and the hyperplane
$x_{uv}=0$ is linearly isomorphic to the cut polytope $\CutP(G/uv)$
where $G/uv$ denotes the contraction of $G$ at the edge $uv$.
We denote by $u$ the node in $G/uv$ representing the edge $uv$ in $G$.
The $uv$-collapsing of an inequality $\vct{a}^\trans\vct{x}\le a_0$ is
an inequality $(\vct{a}')^\trans\vct{x}\le a_0$ defined by
\[
  a'_{ij}=\begin{cases}
    a_{ij}        & \text{if $i,j\ne u$,} \\
    a_{uj}        & \text{if $i=u$, $uj\in E$, $vj\notin E$,} \\
    a_{vj}        & \text{if $i=u$, $uj\notin E$, $vj\in E$,} \\
    a_{uj}+a_{vj} & \text{if $i=u$, $uj,vj\in E$,}
  \end{cases}
\]
for every edge $ij$ of $G/uv$.

The following lemma is given as Lemma~26.4.1~(i) in
\cite{DezLau:cut97}.

\begin{lemma}
  Any collapsing of a valid inequality is valid.
\end{lemma}

\subsubsection{Lifting operations}

The term \emph{lifting} refers to any general operations which derive
an inequality valid for a polyhedron $P$ from an inequality valid for a
polyhedron $P\cap\{\vct{x}\mid x_e=0\}$ for some coordinate
$e$~\cite{NemWol:integer99}.
It is an important way to derive facet inducing inequalities for
combinatorial polyhedra.
In context of the cut polytope, a lifting operation means an operation
which converts an inequality valid for $\CutP(G)$ to an inequality
valid for $\CutP(G')$ where $G$ is obtained by contracting some edges
of $G'$.

Most lifting operations convert an inequality
$\vct{a}^\trans\vct{x}\le a_0$ to an inequality whose appropriate
collapsing is the inequality $\vct{a}^\trans\vct{x}\le a_0$.
Such lifting operations are sometimes called \emph{node splitting}
(see Section~26.5 of \cite{DezLau:cut97}).

The most fundamental lifting operation is zero-lifting.
The following definition and theorem about the zero-lifting of
inequalities for general graphs are due to De~Simone~\cite{Des-ORL90}.

\begin{definition}[Zero-lifting of inequalities]
    \label{def:zero-lifting}
  Let $G=(V,E)$ be a subgraph of $G'=(V',E')$.
  For $\vct{a}\in\RR^E$ and $a_0\in\RR$, the zero-lifting of
  $\vct{a}^\trans\vct{x}\le a_0$ is an inequality
  $(\vct{a}')^\trans\vct{x}\le a_0$ where $\vct{a}'\in\RR^{E'}$ is
  defined by $a'_{uv}=a_{uv}$ for $uv\in E$ and $a'_{uv}=0$ for
  $uv\in E'\setminus E$.
\end{definition}

\begin{theorem} \label{theorem:zero-lifting}
  Let $G=(V,E)$ be a graph with $n$ nodes ($n\ge3$) and $G'=(V',E')$
  be a graph with $n+1$ nodes $V'=V\cup\{w\}$ such that $V$ induces
  $G$ in $G'$.
  Let $(\vct{a}')^\trans\vct{x}\le a_0$ be the zero-lifting of
  $\vct{a}^\trans\vct{x}\le a_0$ and $u$ be a node of $G$.
  Then $(\vct{a}')^\trans\vct{x}\le a_0$ is facet inducing for
  $\CutP(G')$ if the following conditions are met:
  \begin{enumerate}[(i)]
  \item
    $\vct{a}^\trans\vct{x}\le a_0$ is facet inducing for $\CutP(G)$.
  \item
    $\N_{G'}(w)\setminus\{u\}\subseteq\N_G(u)$.
  \item
    The support graph $G(\vct{a})$ of $\vct{a}^\trans\vct{x}\le a_0$
    has at least three nodes.
  \end{enumerate}
\end{theorem}

Theorem~26.5.1 of \cite{DezLau:cut97} is the case of
Theorem~\ref{theorem:zero-lifting} where $G$ and $G'$ are the complete
graphs and $a_0=0$.
The proof of Theorem~26.5.1 of \cite{DezLau:cut97} uses what is called
the lifting lemma (Proposition~2.7 of \cite{DezLau-MP92:facets1} and
Lemma~26.5.3 of \cite{DezLau:cut97}), which has a
wide range of applications.

\begin{lemma}[Lifting lemma] \label{lemma:lifting-1}
  Let $G=(V,E)$ be the complete graph with $n$ nodes $V=\{1,\dots,n\}$
  ($n\ge3$).
  Let $G'=(V',E')$ be the complete graph with $n+1$ nodes
  $V'=V\cup\{n+1\}$.
  Let $\vct{a}\in\RR^E$ and $\vct{a}'\in\RR^{E'}$.
  Suppose that the following assertions hold.
  \begin{enumerate}[(i)]
  \item
    The inequality $\vct{a}^\trans\vct{x}\le0$ is facet inducing for
    $\CUT(G)$ and the inequality $(\vct{a}')^\trans\vct{x}\le0$ is
    valid for $\CUT(G')$.
  \item
    There exist $\abs{E}-1$ subsets $S_j$ of $V\setminus\{1\}$ such
    that the cut vectors $\vct{\delta}_G(S_j)$ are linearly
    independent roots of $\vct{a}^\trans\vct{x}\le0$ and the cut
    vectors $\vct{\delta}_{G'}(S_j)$ are roots of
    $(\vct{a}')^\trans\vct{x}\le0$.
  \item
    There exist $n$ subsets $T_k$ of $V'$ with $1\notin T_k$ and
    $n+1\in T_k$ such that the cut vectors $\vct{\delta}_{G'}(T_k)$
    are roots of $(\vct{a}')^\trans\vct{x}\le0$ and the incidence
    vectors of the sets $T_k$ are linearly independent.
  \end{enumerate}
  Then the inequality $(\vct{a}')^\trans\vct{x}\le0$ is facet inducing
  for $\CUT(G')$.
\end{lemma}

\section{Triangular elimination for general graphs}
  \label{sect:trielim}

\subsection{Definition and validity}

Suppose that we have an inequality $\vct{a}^\trans\vct{x}\le a_0$
which is facet inducing for the cut polytope $\CutP(G)$ of a graph
$G=(V,E)$.
We would like to remove an edge $uv$ from $G$ and instead add some
nodes and edges, converting the inequality
$\vct{a}^\trans\vct{x}\le a_0$ to a facet inducing inequality of
$\CutP(G')$ for the new graph $G'=(V',E')$.

One way to do this is to add a new node $w$ and new edges $uw$ and
$vw$, and add the triangle inequality on $u$, $v$ and $w$ to eliminate
the term $x_{uv}$ from the inequality $\vct{a}^\trans\vct{x}\le a_0$.
For simplicity, we restrict ourselves to the case where $a_0=0$ and
$a_{uv}>0$.
Then the triangle inequality to add is
$-a_{uv}x_{uv}+a_{uv}x_{uw}-a_{uv}x_{vw}\le0$.
This can be seen as a variation of lifting operation since collapsing
the node $w$ to $v$ restores the original inequality, though it
removes an edge from the underlying graph.

\begin{proposition}
  \label{prop:example}
  Let $G=(V,E)$ be a graph and $uv$ be an edge in $G$.
  Let $\vct{a}^\trans\vct{x}\le0$ be a facet inducing inequality of
  $\CUT(G)$ with $a_{uv}>0$.
  Let $w$ be a new node which does not belong to $V$, and let
  $G'=(V',E')$ be a graph with $V'=V\cup\{w\}$ and
  $E'=(E\setminus\{uv\})\cup\{uw,vw\}$.
  Then the inequality
  $\vct{a}^\trans\vct{x}-a_{uv}x_{uv}+a_{uv}x_{uw}-a_{uv}x_{vw}\le0$
  is facet inducing for $\CUT(G')$.
\end{proposition}

Proposition~\ref{prop:example} is a special case of Corollary~2.10~(a)
of \cite{BarMah-MP86}.
We will give a direct proof of Proposition~\ref{prop:example} here
since the proof of Theorem~\ref{theorem:trielim} will follow the same
steps (though with more complicated details).

\begin{proof}
  Let $(\vct{a}')^\trans\vct{x}\le0$ be the new inequality.
  The inequality $(\vct{a}')^\trans\vct{x}\le0$ is valid for
  $\CUT(G'')$, where $G''=(V',E\cup E')$,
  since it is the sum of two inequalities
  $\vct{a}^\trans\vct{x}\le0$ and
  $-a_{uv}x_{uv}+a_{uv}x_{uw}-a_{uv}x_{vw}\le0$
  both of which are valid for $\CUT(G'')$.
  The inequality $(\vct{a}')^\trans\vct{x}\le0$ is also valid for
  $\CUT(G')$ since it consists of terms corresponding to edges of
  $G'$, which is a subgraph of $G''$.

  Since $\vct{a}^\trans\vct{x}\le0$ is facet inducing for
  $\CUT(G)$, there exist $\abs{E}-1$ subsets $S_1,\dots,S_{\abs{E}-1}$
  of $V\setminus\{v\}$ such that the $\abs{E}-1$ cut vectors
  $\vct{\delta}_G(S_j)$ are linearly independent roots of
  $\vct{a}^\trans\vct{x}\le0$.

  If we collapse the node $w$ to the node $v$ in
  $(\vct{a}')^\trans\vct{x}\le0$, we obtain the inequality
  $\vct{a}^\trans\vct{x}\le0$.
  This implies that $\vct{\delta}_{G'}(S_j)$ are linearly independent
  roots of $(\vct{a}')^\trans\vct{x}\le0$.
  The $\abs{E}-1$ cut vectors $\vct{\delta}_{G'}(S_j)$ satisfy an
  equation $x_{vw}=0$.
  On the other hand, a cut vector $\vct{\delta}_{G'}(\{w\})$ is a root
  of $(\vct{a}')^\trans\vct{x}\le0$ with $x_{vw}=1\ne0$.
  Therefore, the $\abs{E}=\abs{E'}-1$ roots $\vct{\delta}_{G'}(S_j)$
  and $\vct{\delta}_{G'}(\{w\})$ of $(\vct{a}')^\trans\vct{x}\le0$ are
  linearly independent.
  This implies that $(\vct{a}')^\trans\vct{x}\le0$ is facet inducing
  for $\CUT(G')$.
\end{proof}

A special case of Theorem~\ref{theorem:cycle}~(ii)
where the graph $G$ is identical to the cycle $C$ may be proved by
using Proposition~\ref{prop:example} repeatedly as follows.

\begin{corollary}
  The cycle inequality~(\ref{eq:cycle}) is facet inducing for
  $\CUT(C)$.
\end{corollary}

\begin{proof}
  The proof is by induction on the length $n$ of the cycle $C$.
  If $n=3$, then the inequality~(\ref{eq:cycle}) is the triangle
  inequality and facet inducing for $\CUT(C)$.
  In case of $n>3$, we let $v$ be the node in $C$ adjacent of $w$
  other than $u$ and apply Proposition~\ref{prop:example} with
  $G=C/vw$, $G'=C$, and $\vct{a}^\trans\vct{x}\le0$ is the cycle
  inequality in $C/vw$, which is facet inducing for $C/vw$ by the
  induction hypothesis.
\end{proof}

One question arises here: can we add more edges to $G'$ keeping the
property that the new inequality is facet inducing for $\CutP(G')$?
We will answer this question affirmatively by
Theorem~\ref{theorem:trielim}.
The main ingredient of the proof is the notion of
triangular elimination, which generalizes the operation
described in Proposition~\ref{prop:example}.

In what follows, we use the following notation and terms.
Let $\Tri(u,v;w)=x_{uv}-x_{uw}-x_{vw}$ and
$\Tri(u,v,w)=x_{uv}+x_{uw}+x_{vw}-2$ for any three nodes $u,v,w$ in
the graph in question.
The notation $\Tri\{u,v,w\}$ ambiguously denotes one of the four
triangular forms $\Tri(u,v;w)$, $\Tri(w,v;u)$, $\Tri(u,w;v)$ or
$\Tri(u,v,w)$.
The \emph{support graph} of a vector $\vct{a}\in\RR^E$ is a subgraph
$G(\vct{a})=(V(\vct{a}),E(\vct{a}))$ of $G$ whose edges are all edges
$e$ in $G$ with $a_e\ne0$ and nodes are all the endpoints of the edges
in $E(\vct{a})$.
For a vector $\vct{a}\in\RR^E$, a scalar $a_0\in\RR$ and a subset
$F\subseteq E$, we say the inequality $\vct{a}^\trans\vct{x}\le a_0$
is \emph{completely supported} by $F$ when $E(\vct{a})$ is included in
$F$.

\begin{definition}[Triangular elimination for graphs]
    \label{def:trielim-graph}
  Let $G=(V,E)$ be a graph, $t$ be an integer, and let
  $F=\{u_iv_i\mid i=1,\dots,t\}$ be any subset of $E$.
  The graph $G'=(V',E')$ is a \emph{triangular elimination} of $G$
  (with respect to $F$) if $V'=V\cup\{w_1,\dots,w_t\}$,
  $E'\supseteq\{w_iu_i,w_iv_i\mid i=1,\dots,t\}$, and
  $E'\cap E=E\setminus F$.
  Here $w_1,\dots,w_t$ are distinct nodes not in $V$.
  Node $w_i$ of $G'$ is said to be \emph{associated} with edge $u_iv_i$
  of $G$.
\end{definition}

\begin{definition}[Triangular elimination for inequalities]
    \label{def:trielim-ineq}
  Let $G'=(V',E')$ be a triangular elimination of $G=(V,E)$, and
  suppose we are given $\vct{a}\in\RR^E$, $a_0\in\RR$,
  $\vct{a}'\in\RR^{E'}$, $a'_0\in\RR$.
  Then inequality $(\vct{a}')^\trans\vct{x}\le a'_0$ is a
  \emph{triangular elimination} of $\vct{a}^\trans\vct{x}\le a_0$ if
  for some choices of triangular forms $\Tri_i\{u_i,v_i,w_i\}$,
  $i=1,\dots,t$, we have
  \[
    (\vct{a}')^\trans\vct{x}-a'_0=\vct{a}^\trans\vct{x}-a_0
    +\sum_{i=1}^t \abs{a_{u_iv_i}}\Tri_i\{u_i,v_i,w_i\}.
  \]
\end{definition}

The operation in Proposition~\ref{prop:example} is the case where
$t=1$, $u_1=u$, $v_1=v$, $w_1=w$, $\Tri_1\{u,v,w\}=\Tri(u,w;v)$, and
$w$ has no neighbours other than $u$ and $v$ in $G'$.

\begin{figure}
  \centering
  \includegraphics{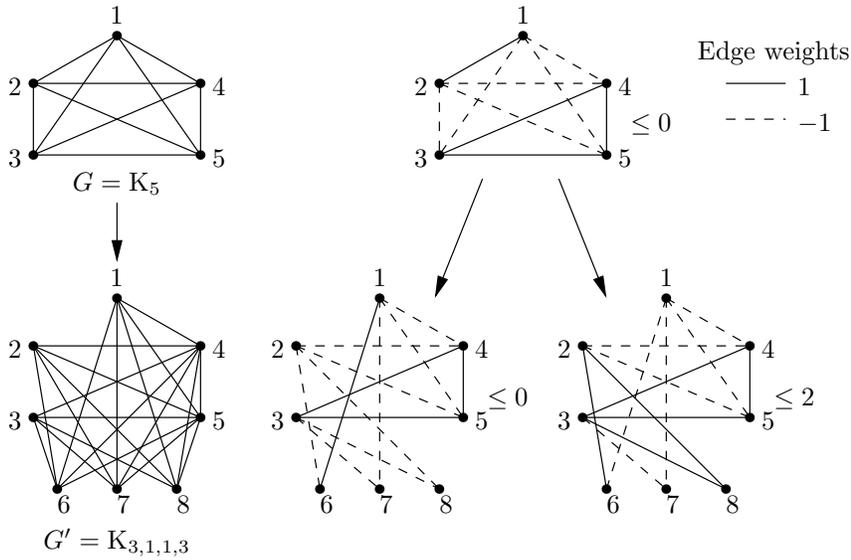}
  \caption{A graph $G'$ and two inequalities obtained as the triangular
    eliminations of $G=\K_5$ and its facet inducing pentagonal
    inequality.} \label{fig:trielim-ex}
\end{figure}

\begin{example} \label{example:trielim-ex}
  To understand how triangular elimination typically works, let us
  consider another, more explicit, example (see
  Figure~\ref{fig:trielim-ex}).
  Let $G=(V,E)$ be the complete graph $\K_5$ with nodes labelled as
  $1,\dots,5$ and $G'=(V',E')$ be the complete 4-partite graph
  $\K_{3,1,1,3}$ with vertices partitioned as $\{1,2,3\}$, $\{4\}$,
  $\{5\}$ and $\{6,7,8\}$.
  Then $G'$ is a triangular elimination of $G$ with respect to
  $F=E\setminus E'=\{12,13,23\}$, where $u_1v_1=12$, $u_2v_2=13$,
  $u_3v_3=23$ and $w_i=5+i$ for $i=1,2,3$.
  Let $\vct{a}^\trans\vct{x}\le a_0$ be the pentagonal inequality
  $x_{12}+x_{34}+x_{35}+x_{45}-\sum_{1\le u\le2,3\le v\le5}x_{uv}
   \le0$.
  Then one of the triangular eliminations
  $(\vct{a}')^\trans\vct{x}\le a'_0$ of $\vct{a}^\trans\vct{x}\le a_0$
  is obtained by using triangular forms
  $\Tri\{1,2,6\}=\Tri(1,6;2)$, $\Tri\{1,3,7\}=\Tri(1,3;7)$ and
  $\Tri\{2,3,8\}=\Tri(2,3;8)$, and it is
  $-\sum_{1\le u\le2,4\le v\le5}x_{uv}+x_{34}+x_{35}+x_{45}
   +x_{16}-x_{26}-x_{17}-x_{37}-x_{28}-x_{38}\le0$.
  Another triangular elimination $(\vct{a}'')^\trans\vct{x}\le a''_0$
  of $\vct{a}^\trans\vct{x}\le a_0$ is obtained by using
  $\Tri\{1,2,6\}=\Tri(2,6;1)$, $\Tri\{1,3,7\}=\Tri(1,3;7)$ and
  $\Tri\{2,3,8\}=\Tri(2,3,8)$, and it is
  $-\sum_{1\le u\le2,4\le v\le5}x_{uv}+x_{34}+x_{35}+x_{45}
   -x_{16}+x_{26}-x_{17}-x_{37}+x_{28}+x_{38}\le2$.
\end{example}

Triangular elimination of an inequality can also be seen as a specific
combination of zero-lifting operation and Fourier-Motzkin elimination.
Let $\vct{a}^\trans\vct{x}\le a_0$ be a valid inequality of
$\CutP(G)$.
Consider a graph $G''=(V',E\cup E')$ and the zero-lifting of
$\vct{a}^\trans\vct{x}\le a_0$ to $\CutP(G'')$.
Then apply Fourier-Motzkin elimination to project out the variables
$x_{uv}$ for $uv\in F$, adding triangle inequalities to
$\vct{a}^\trans\vct{x}\le a_0$.
This gives a triangular elimination of $\vct{a}^\trans\vct{x}\le a_0$.

The ``if'' part of the following theorem is straightforward from the
definitions.

\begin{theorem} \label{theorem:valid}
  Let $G'=(V',E')$ be a triangular elimination of $G=(V,E)$, and let
  $(\vct{a}')^\trans\vct{x}\le a'_0$ be a triangular elimination of
  $\vct{a}^\trans\vct{x}\le a_0$.
  Then $(\vct{a}')^\trans\vct{x}\le a'_0$ is valid for $\CutP(G')$ if
  and only if $\vct{a}^\trans\vct{x}\le a_0$ is valid for $\CutP(G)$.
\end{theorem}

\begin{proof}[Proof of the ``if'' part of Theorem~\ref{theorem:valid}]
  Let $G''=(V',E\cup E')$.
  The inequality $(\vct{a}')^\trans\vct{x}\le a'_0$ is
  valid for $\CutP(G'')$ since it is a sum of an
  inequality $\vct{a}^\trans\vct{x}\le a_0$ and triangle
  inequalities all of which are valid for $\CutP(G'')$.
  The inequality $(\vct{a}')^\trans\vct{x}\le a'_0$ is also valid for
  $\CutP(G')$ since it consists of terms corresponding to edges of
  $G'$, which is a subgraph of $G''$.
\end{proof}

We prove the ``only if'' part in Section~\ref{subsect:switching}.

\subsection{Switching of the triangular elimination}
  \label{subsect:switching}

As is shown in Example~\ref{example:trielim-ex},
Definition~\ref{def:trielim-ineq} allows several choices of
$\Tri_i\{u_i,v_i,w_i\}$, and different choices give apparently
different inequalities.
This may complicate handling of triangular elimination.
It turns out that if we deal with equivalence classes of inequalities
under switching equivalence instead of the inequalities themselves,
triangular elimination is easier to handle.

\begin{proposition} \label{prop:switching}
  Let $G'$ be a triangular elimination of $G$.
  Let $(\vct{a}')^\trans\vct{x}\le a'_0$ be a triangular elimination
  of $\vct{a}^\trans\vct{x}\le a_0$ and
  $(\vct{b}')^\trans\vct{x}\le b'_0$ be a triangular elimination of
  $\vct{b}^\trans\vct{x}\le b_0$ such that the association of nodes in
  $G'$ with edges in $G$ are the same in both triangular eliminations.
  If $\vct{a}^\trans\vct{x}\le a_0$ is switching equivalent to
  $\vct{b}^\trans\vct{x}\le b_0$, then
  $(\vct{a}')^\trans\vct{x}\le a'_0$ is switching equivalent to
  $(\vct{b}')^\trans\vct{x}\le b'_0$.
\end{proposition}

\begin{proof}
  Let
  \[
    ((\vct{a}')^\trans\vct{x}-a'_0)-(\vct{a}^\trans\vct{x}-a_0)
    =\sum_{1\le i\le t}\abs{a_{u_iv_i}}\Tri_i\{u_i,v_i,w_i\}.
  \]
  First we prove the proposition when $\vct{a}=\vct{b}$ and $a_0=b_0$.
  In this case, let
  \[
    ((\vct{b}')^\trans\vct{x}-b'_0)-(\vct{a}^\trans\vct{x}-a_0)
    =\sum_{1\le i\le t}\abs{a_{u_iv_i}}\tilde{\Tri}_i\{u_i,v_i,w_i\}.
  \]
  For $i=1,\dots,t$, if $a_{u_iv_i}\ne0$, then $\Tri_i\{u_i,v_i,w_i\}$
  and $\tilde{\Tri}_i\{u_i,v_i,w_i\}$ are either identical or the
  $\{w_i\}$-switching of each other, by comparing their
  $x_{u_iv_i}$-coefficient.
  Let $S$ be the set of $w_i$ for $i$ such that
  $\Tri_i\{u_i,v_i,w_i\}$ and $\tilde{\Tri}_i\{u_i,v_i,w_i\}$ are the
  $\{w_i\}$-switching of each other.
  Then two inequalities $(\vct{a}')^\trans\vct{x}\le a'_0$ and
  $(\vct{b}')^\trans\vct{x}\le b'_0$ are the $S$-switching of each
  other.

  Next we prove the general case.
  Let $S$ be a subset of $V$ such that $\vct{a}^\trans\vct{x}\le a_0$
  and $\vct{b}^\trans\vct{x}\le b_0$ are the $S$-switching of each
  other.
  Let $(\vct{b}'')^\trans\vct{x}\le b''_0$ be the $S$-switching of
  $(\vct{a}')^\trans\vct{x}\le a'_0$.
  Then
  $((\vct{b}'')^\trans\vct{x}-b''_0)-(\vct{b}^\trans\vct{x}-b_0)\le0$
  is the $S$-switching of the inequality
  $((\vct{a}')^\trans\vct{x}-a'_0)-(\vct{a}^\trans\vct{x}-a_0)\le0$
  and therefore the $S$-switching of
  $\sum_{1\le i\le t}\abs{a_{u_iv_i}}\Tri_i\{u_i,v_i,w_i\}\le0$.
  This means that $(\vct{b}'')^\trans\vct{x}\le b''_0$ as well as
  $(\vct{b}')^\trans\vct{x}\le b'_0$ is a triangular elimination of
  $\vct{b}^\trans\vct{x}\le b_0$.
  By the case we already proved, $(\vct{b}')^\trans\vct{x}\le b'_0$
  and $(\vct{b}'')^\trans\vct{x}\le b''_0$ are switching equivalent.
  Therefore, $(\vct{a}')^\trans\vct{x}\le a'_0$ and
  $(\vct{b}')^\trans\vct{x}\le b'_0$ are switching equivalent.
\end{proof}

\begin{example}[continued from Example~\ref{example:trielim-ex}]
  Both inequalities $(\vct{a}')^\trans\vct{x}\le a'_0$ and
  $(\vct{a}'')^\trans\vct{x}\le a''_0$ described in
  Example~\ref{example:trielim-ex}) are the triangular eliminations of
  $\vct{a}^\trans\vct{x}\le a_0$.
  By Proposition~\ref{prop:switching},
  $(\vct{a}')^\trans\vct{x}\le a'_0$ and
  $(\vct{a}'')^\trans\vct{x}\le a''_0$ are switching equivalent.
  In fact, they are the $\{6,8\}$-switching equivalent of each other.
\end{example}

Proposition~\ref{prop:switching} essentially states that triangular
elimination is well-defined as an operation acting on
switching-equivalence classes of inequalities.
By Proposition~\ref{prop:switching}, we can freely replace
$\vct{a}^\trans\vct{x}\le a_0$ with its switching and we do not need
to care the choice of $\Tri_i$ when we
are interested in switching-invariant properties of the inequalities
obtained by triangular elimination such as whether it is valid or not,
facet inducing or not, and so on.
Any properties of inequalities that we deal with in the rest of the
paper are switching-invariant.

By using Proposition~\ref{prop:switching}, we can now
complete the proof of
Theorem~\ref{theorem:valid}.

\begin{proof}[Proof of the ``only if'' part of Theorem~\ref{theorem:valid}]
  Suppose that $(\vct{a}')^\trans\vct{x}\le a'_0$ is valid for
  $\CutP(G')$.
  By Proposition~\ref{prop:switching}, we can assume without loss of
  generality that for $1\le i\le t$,
  $\Tri_i\{u_i,v_i,w_i\}=\Tri(u_i,v_i;w_i)$ (if $a_{u_iv_i}\le0$) or
  $\Tri_i\{u_i,v_i,w_i\}=\Tri(u_i,w_i;v_i)$ (if $a_{u_iv_i}\ge0$).
  Then the inequality $\vct{a}^\trans\vct{x}\le a_0$ is
  obtained from $(\vct{a}')^\trans\vct{x}\le a'_0$ by collapsing the
  node $w_i$ to $v_i$ for every $1\le i\le t$.
  This means that the inequality $\vct{a}^\trans\vct{x}\le a_0$ is
  also valid.
\end{proof}

\subsection{Facets and triangular elimination}

We state and prove a sufficient condition for triangular elimination
to be facet preserving.
Note the similarity to the conditions in
Theorem~\ref{theorem:zero-lifting}.

\begin{theorem} \label{theorem:trielim}
  Let $G'=(V',E')$ be a triangular elimination of $G=(V,E)$, and let
  $(\vct{a}')^\trans\vct{x}\le a'_0$ be a triangular elimination of
  $\vct{a}^\trans\vct{x}\le a_0$.
  Then $(\vct{a}')^\trans\vct{x}\le a'_0$ is facet inducing for
  $\CutP(G')$ if the following conditions apply:
  \begin{enumerate}[(i)]
  \item
    The inequality $\vct{a}^\trans\vct{x}\le a_0$ is facet inducing
    for $\CutP(G)$.
  \item
    For $i=1,\dots,t$,
    $\N_{G'}(w_i)\setminus\{u_i,v_i\}\subseteq\N_G(u_i)\cap\N_G(v_i)$.
  \item
    For $i=1,\dots,t$, the inequality $\vct{a}^\trans\vct{x}\le a_0$
    is not completely supported by the edge set
    $\{u_il,v_il\mid l\in\N_{G'}(w_i)\}$.
  \end{enumerate}
\end{theorem}

Note that condition~(ii) implies that the set $\{w_i\mid i=1,\dots,t\}$
is an independent set in $G'$.

\begin{example}[continued from Example~\ref{example:trielim-ex}]
  The inequality $(\vct{a}')^\trans\vct{x}\le a'_0$ described in
  Example~\ref{example:trielim-ex} is facet inducing for $\CutP(G')$,
  since the graphs $G,G'$ and the inequalities
  $\vct{a}^\trans\vct{x}\le a_0,(\vct{a}')^\trans\vct{x}\le a'_0$
  satisfy the conditions in Theorem~\ref{theorem:trielim}.
\end{example}

We prove Theorem~\ref{theorem:trielim} in a similar way to the proof of
the zero-lifting theorem, Theorem~26.5.1 of \cite{DezLau:cut97}, as
follows.
We first introduce a variation of the lifting lemma,
Lemma~\ref{lemma:lifting-1}, adapted to graphs other than the
complete graphs and contraction of multiple edges.
Then to prove one of the preconditions of the lemma, we use a lemma
from \cite{DezLau:cut97}.

First we introduce the variation of the lifting lemma.

\begin{lemma} \label{lemma:lifting}
  Let $G'=(V',E')$ be a graph and $H=(V',F)$ be a forest in $G'$ with
  $t$ edges $F=\{v_1w_1,\dots,v_tw_t\}\subseteq E'$.
  Let $G=(V,E)$ be the graph obtained from $G'$ by contracting the
  edges in $H$.
  Let $U_i=\N_{G'}(v_i)\cap\N_{G'}(w_i)$.
  We require that $\abs{E'}=\abs{E}+\abs{U_1}+\dots+\abs{U_t}+t$.
  Let $\vct{a}\in\RR^E$ and $\vct{a}'\in\RR^{E'}$.
  Suppose that the following assertions hold.
  \begin{enumerate}[(i)]
  \item
    The inequality $\vct{a}^\trans\vct{x}\le0$ is facet inducing for
    $\CUT(G)$ and the inequality $(\vct{a}')^\trans\vct{x}\le0$ is
    valid for $\CUT(G')$.
  \item
    There exist $\abs{E}-1$ subsets $\tilde{S}_j$ of $V$ such that the
    cut vectors $\vct{\delta}_G(\tilde{S}_j)$ are linearly independent
    roots of $\vct{a}^\trans\vct{x}\le0$ and the cut vectors
    $\vct{\delta}_{G'}(S_j)$ are roots of
    $(\vct{a}')^\trans\vct{x}\le0$, where
    $S_j=\tilde{S}_j\cup\{w_i\mid v_i\in\tilde{S}_j\}$.
  \item
    For $1\le i\le t$, there exist $\abs{U_i}+1$ subsets $T_{ik}$ of
    $V'$ with $v_i\notin T_{ik}$, $w_i\in T_{ik}$ and
    $\delta_{v_lw_l}(T_{ik})=0$ for $1\le l\le t,\;l\ne i$ such that
    the cut vectors $\vct{\delta}_{G'}(T_{ik})$ are roots of
    $(\vct{a}')^\trans\vct{x}\le0$ and the incidence vectors of the
    sets $T_{ik}\cap(U_i\cup\{w_i\})$ are linearly independent.
  \end{enumerate}
  Then the inequality $(\vct{a}')^\trans\vct{x}\le0$ is facet inducing
  for $\CUT(G')$.
\end{lemma}

Note that Lemma~\ref{lemma:lifting-1} is a
special case of this lemma with $t=1$, $G=\K_n$, $G'=\K_{n+1}$,
$v_1=1$ and $w_1=n+1$.
The proof is similar to the latter half of the proof of Theorem~26.5.1
of \cite{DezLau:cut97}, though our proof is a little more complicated
because we cannot use a correlation cone instead of the cut cone
$\CUT(G')$.

\begin{remark}
  The same remark on node splitting as that given below Lemma~26.5.3
  of \cite{DezLau:cut97} applies for Lemma~\ref{lemma:lifting}.
  That is, if the inequality $\vct{a}^\trans\vct{x}\le0$ comes from
  $(\vct{a}')^\trans\vct{x}\le0$ by collapsing the nodes $w_i$ to the
  corresponding nodes $v_i$, then the assertion~(ii) is implied by the
  assertion~(i).
\end{remark}

\begin{proof}
  Note that $\abs{E'}-1=(\abs{E}-1)+(\abs{U_1}+1)+\dots+(\abs{U_t}+1)$.
  We show that $\abs{E'}-1$ cut vectors $\vct{\delta}_{G'}(S_j)$ and
  $\vct{\delta}_{G'}(T_{ik})$ are linearly independent.
  Let us consider the $\abs{E'}\times(\abs{E'}-1)$ matrix $M$, whose
  columns are these $\abs{E'}-1$ cut vectors.
  We prove that $M$ has full column rank.
  The rows of $M$ are indexed by the edges in $E'$, which can be
  grouped as $E'=I\cup\bigcup_{1\le i\le t}(J_i\cup K_i\cup L_i)$:
  \begin{itemize}
  \item $I$ consists the edges in $E'$ which do not belong to any of
    the following groups,
  \item $J_i=\{v_i u\mid u\in U_i\}$,
  \item $K_i=\{w_i u\mid u\in U_i\}$,
  \item $L_i=\{v_iw_i\}$.
  \end{itemize}
  Note that some edges in $E'$ may belong to more than one set.
  In that case, we consider that $M$ contains the corresponding rows
  twice.
  We can do so since this does not change the rank of $M$.
  Then the matrix $M$ is of the form:
  \[
    \def\vi{\vct{1}}
    M=\begin{blockarray}{ssssss}
             & (S_j)  &  (T_{1k})  &  (T_{2k})  & \cdots &  (T_{tk})  \\
      \begin{block}{s(ccccc)}
       (I)   &  X_0   &    X_1     &    X_2     & \cdots &    X_t     \\
      (J_1)  & Y_{01} &   Y_{11}   &   Y_{21}   & \cdots &   Y_{t1}   \\
      (J_2)  & Y_{02} &   Y_{12}   &   Y_{22}   & \cdots &   Y_{t2}   \\
      \vdots & \vdots &   \vdots   &   \vdots   & \ddots &   \vdots   \\
      (J_t)  & Y_{0t} &   Y_{1t}   &   Y_{2t}   & \cdots &   Y_{tt}   \\
      (K_1)  & Y_{01} & \vi-Y_{11} &   Y_{21}   & \cdots &   Y_{t1}   \\
      (L_1)  &   0    &    \vi     &     0      & \cdots &     0      \\
      (K_2)  & Y_{02} &   Y_{12}   & \vi-Y_{22} & \cdots &   Y_{t2}   \\
      (L_2)  &   0    &     0      &    \vi     & \cdots &     0      \\
      \vdots & \vdots &   \vdots   &   \vdots   & \ddots &   \vdots   \\
      (K_t)  & Y_{0t} &   Y_{1t}   &   Y_{2t}   & \cdots & \vi-Y_{tt} \\
      (L_t)  &   0    &     0      &     0      & \cdots &    \vi     \\
      \end{block}
    \end{blockarray},
  \]
  where $\vct{1}$ denotes the all-ones matrix.
  To prove that $M$ has full column rank, we transform $M$ by
  reversible linear operations on its row vectors as follows: subtract
  the rows corresponding to the edge $v_i u$ in $J_i$ from the rows
  corresponding to the edge $w_i u$ in $K_i$, subtract the row $L_i$
  from each row in $K_i$, and divide the rows in $K_i$ by $-2$.
  Then we obtain:
  \[
    \def\vi{\vct{1}}
    M'=\begin{blockarray}{ssssss}
             & (S_j)  & (T_{1k}) & (T_{2k}) & \cdots & (T_{tk}) \\
      \begin{block}{s(ccccc)}
       (I)   &  X_0   &   X_1    &   X_2    & \cdots &   X_t    \\
      (J_1)  & Y_{01} &  Y_{11}  &  Y_{21}  & \cdots &  Y_{t1}  \\
      (J_2)  & Y_{02} &  Y_{12}  &  Y_{22}  & \cdots &  Y_{t2}  \\
      \vdots & \vdots &  \vdots  &  \vdots  & \ddots &  \vdots  \\
      (J_t)  & Y_{0t} &  Y_{1t}  &  Y_{2t}  & \cdots &  Y_{tt}  \\
      (K_1)  &   0    &  Y_{11}  &    0     & \cdots &    0     \\
      (L_1)  &   0    &   \vi    &    0     & \cdots &    0     \\
      (K_2)  &   0    &    0     &  Y_{22}  & \cdots &    0     \\
      (L_2)  &   0    &    0     &   \vi    & \cdots &    0     \\
      \vdots & \vdots &  \vdots  &  \vdots  & \ddots &  \vdots  \\
      (K_t)  &   0    &    0     &    0     & \cdots &  Y_{tt}  \\
      (L_t)  &   0    &    0     &    0     & \cdots &   \vi    \\
      \end{block}
    \end{blockarray}.
  \]
  The leftmost $\abs{E}-1$ columns of $M'$ have full column rank by
  assertion~(ii).
  The $k$th column of the square submatrix
  $\bigl(\begin{smallmatrix}Y_{ii} \\ \vct{1}\end{smallmatrix}\bigr)$
  of order $\abs{U_i}+1$ of $M'$ is the incidence vector of the set
  $T_{ik}\cap(U_i\cup\{w_i\})$.
  This implies that
  $\bigl(\begin{smallmatrix}Y_{ii} \\ \vct{1}\end{smallmatrix}\bigr)$
  has full rank by assertion~(iii).
  Therefore, $M'$ (and also $M$) has full column rank.
\end{proof}

The other lemma which we use is Lemma~26.5.2~(ii) of
\cite{DezLau:cut97}.
In \cite{DezLau:cut97} the graph $G$ is restricted to the complete
graph, but this restriction is not relevant.

\begin{lemma}[\cite{DezLau:cut97}]
  \label{lemma:projection}
  Let $G=(V,E)$ be a graph and $\vct{a}^\trans\vct{x}\le0$ be an
  inequality inducing a facet $F$ of $\CUT(G)$.
  Let $D$ be a subset of $E$ and $F'$ be the projection of
  $F\subseteq\RR^E$ to $\RR^D$.
  If there exists an edge $e\in E\setminus D$ with $a_e\ne0$, then
  $F'$ is full-dimensional.
  Otherwise, the dimension of $F'$ is $\abs{D}-1$.
\end{lemma}

Using Lemmas~\ref{lemma:lifting} and \ref{lemma:projection}, we
prove Theorem~\ref{theorem:trielim}.

\begin{proof}[Proof of Theorem~\ref{theorem:trielim}]
  By Propositions~\ref{prop:cone-switching} and \ref{prop:switching},
  we can assume without loss of
  generality that $a_0$ is equal to zero and that for $1\le i\le t$,
  $\Tri_i\{u_i,v_i,w_i\}$ is equal to
  $\Tri(u_i,v_i;w_i)$ (if $a_{u_iv_i}\le0$) or
  $\Tri(u_i,w_i;v_i)$ (if $a_{u_iv_i}\ge0$).
  These assumptions imply that neither $\vct{a}^\trans\vct{x}\le a_0$
  nor $(\vct{a}')^\trans\vct{x}\le a'_0$ has a nonzero constant term.
  By Proposition~\ref{prop:cone-homogeneous}, this means that
  $\vct{a}^\trans\vct{x}\le0$ is facet inducing for
  the cut cone $\CUT(G)$, and that $(\vct{a}')^\trans\vct{x}\le0$
  is valid for the cut cone $\CUT(G')$.

  We prove that $(\vct{a}')^\trans\vct{x}\le0$ induces a facet of
  $\CUT(G')$ by using Lemma~\ref{lemma:lifting}.

  The assertion~(i) holds by Theorem~\ref{theorem:valid}.
  The assertion~(ii) holds since the inequality
  $\vct{a}^\trans\vct{x}\le0$ comes from the inequality
  $(\vct{a}')^\trans\vct{x}\le0$ by collapsing the nodes $w_i$ to the
  corresponding nodes $v_i$.
  All we need to do is to check the assertion~(iii).

  Let $1\le i\le t$ and $m=\abs{U_i}$.
  We define $D_i=\{u_iv_i\}\cup\{u_il,v_il\in E\mid l\in U_i\}$.
  We construct $m+1$ subsets $T_{ik}$ of $V'$ satisfying the
  assertion~(iii).

  Let $F$ be the facet of $\CUT(G)$ induced by
  $\vct{a}^\trans\vct{x}\le0$.
  By Lemma~\ref{lemma:projection} with $D=D_i$, the projection $F_i$
  of $F$ to $\RR^{D_i}$ is full-dimensional.
  This means that we have $\abs{D_i}=2m+1$ subsets $\tilde{T}_{ik}$ of
  $V$ with $v_i\notin\tilde{T}_{ik}$ such that
  $\vct{\delta}_G(\tilde{T}_{ik})$ are roots of
  $\vct{a}^\trans\vct{x}\le0$ and the $2m+1$ cut vectors
  $\vct{\delta}_{D_i}(\tilde{T}_{ik})$ are linearly independent.

  We show that $u_i$ belongs to exactly $m+1$ out of the $2m+1$ sets
  $\tilde{T}_{ik}$.
  To show this by contradiction, first suppose that $u_i$ belongs to
  at most $m$ of them.
  This means that at least $m+1$ of them does not contain $u_i$ and
  that the intersection
  $F_i\cap\{\vct{x}\in\RR^{D_i}\mid x_{u_iv_i}=0\}$ has a dimension at
  least $m+1$.
  However, this intersection is contained in the intersection
  $\CUT(D_i)\cap\{\vct{x}\in\RR^{D_i}\mid x_{u_iv_i}=0\}$,%
  \footnote{Here we denote by $\CUT(D_i)$ the cut cone of the graph
    with edges $D_i$ and nodes $V$ or any subset of $V$ that contains
    $U_i\cup\{u_i,v_i\}$.
    We justify this slight abuse of the notation by the fact that
    adding isolated nodes to a graph does not change the cut cone.}
  whose dimension is $m$, a contradiction.
  Thus $u_i$ belongs to at least $m+1$ out of the $2m+1$ subsets.
  On the other hand, suppose that $u_i$ belongs to at least $m+2$ of
  them.
  If $u_i\in\tilde{T}_{ik}$, then $\vct{\delta}_{D_i}(\tilde{T}_{ik})$
  satisfies equations $x_{u_iv_i}=x_{u_il}+x_{v_il}$ for all
  $l\in U_i$.
  This implies that the $(m+1)$-dimensional subspace of $\RR^{D_i}$
  defined by $x_{u_iv_i}=x_{u_il}+x_{v_il}$ for $l\in U_i$ contains
  $m+2$ linearly independent vectors, a contradiction.
  Thus $u_i$ belongs to exactly $m+1$ out of the $2m+1$ sets
  $\tilde{T}_{ik}$.
  As a result, we can assume
  $u_i\in\tilde{T}_{i1},\dots,\tilde{T}_{i,m+1}$ and
  $u_i\notin\tilde{T}_{i,m+2},\dots,\tilde{T}_{i,2m+1}$ without loss
  of generality.
  We define $m+1$ subsets $T_{ik}$ of $V'$ as follows.
  If $a_{u_iv_i}\le0$, then let
  $T_{ik}=\tilde{T}_{ik}\cup
          \{w_l\mid v_l\in\tilde{T}_{ik}\}\cup\{w_i\}$
  for $1\le k\le m+1$.
  Otherwise, let
  $T_{ik}=\tilde{T}_{i,m+1+k}\cup
          \{w_l\mid v_l\in\tilde{T}_{i,m+1+k}\}\cup\{w_i\}$
  for $1\le k\le m$ and $T_{i,m+1}=\{w_i\}$.

  Now we prove that the incidence vectors of $m+1$ sets
  $T_{ik}\cap(U_i\cup\{w_i\})$ are linearly independent.
  Let $M$ be the $(2m+1)\times(m+1)$ matrix whose $k$th column vector
  is the cut vector $\vct{\delta}_{D_i}(T_{ik})$,
  and $M'$ be the square matrix of order $m+1$ whose $k$th column
  vector is the incidence vector of $T_{ik}\cap(U_i\cup\{w_i\})$.
  We prove that $M'$ is nonsingular.
  The matrix $M'$ is of the form
  $M'=\bigl(\begin{smallmatrix}X \\ \vct{1}\end{smallmatrix}\bigr)$,
  where the bottommost row corresponds to the node $w_i$.
  The rows of $M$ correspond to the edges in $D_i$ which are grouped
  as $D_i=J\cup K\cup L$: $J=\{u_il\mid l\in U_i\}$,
  $K=\{v_il\mid l\in U_i\}$ and $L=\{u_iv_i\}$.
  If $a_{u_iv_i}\le0$, then the matrix $M$ is given by:
  \[
    M=\begin{blockarray}{s(c)}
      (J) & \vct{1}-X \\
      (K) &     X     \\
      (L) &  \vct{1}
    \end{blockarray},
  \]
  and it has full column rank by assumption.
  Without decreasing its rank, we can transform $M$ to $M'$
  by reversible linear operations on rows and removing
  all-zero rows.
  This means $M'$ is nonsingular.
  Similarly, if $a_{u_iv_i}>0$, then the matrix $M$ is given by:
  \[
    M=\begin{blockarray}{s(cc)}
      (J) &     X     \\
      (K) &     X     \\
      (L) &  \vct{0}
    \end{blockarray}.
  \]
  Its leftmost $m$ columns are linearly independent and its rightmost
  column is the all-zero vector by assumption.
  By a similar argument as above, the leftmost $m$ columns of $X$ are
  linearly independent.
  This implies the $m+1$ column vectors of $X$ are affinely
  independent, or equivalently the matrix $M'$ is nonsingular.
  This means the assertion~(iii) is satisfied.
\end{proof}

To make Theorem~\ref{theorem:trielim} easier to use, we show that
condition~(iii) in Theorem~\ref{theorem:trielim} holds for any facet
inducing inequalities except for the triangle inequality and the
inequality of the forms $x_e\ge0$ and $x_e\le1$.

\begin{proposition} \label{prop:not-triangle}
  Let $G=(V,E)$ be a graph and $\vct{a}^\trans\vct{x}\le a_0$ be a
  facet inducing inequality of $\CutP(G)$.
  Let $u_1v_1,\dots,u_tv_t$ be $t$ distinct edges of $G$.
  If the support graph of the vector $\vct{a}$ has more than three
  nodes, then the inequality $\vct{a}^\trans\vct{x}\le a_0$ is not
  completely supported by the edge set
  $\{u_il,v_il\mid l\in\N_G(u_i)\cap\N_G(v_i)\}$ for any
  $i=1,\dots,t$.
\end{proposition}

To prove Proposition~\ref{prop:not-triangle}, we need the
following lemma.

\begin{lemma} \label{lemma:degree2}
  Let $G=(V,E)$ be a graph, $\vct{a}\in\RR^E$ be a vector, and
  $a_0\in\RR$ be a scalar.
  Suppose the following assumptions hold.
  \begin{enumerate}[(i)]
  \item
    $G$ contains a triangle on nodes $l,u,v$ as a subgraph.
  \item
    At least one of $a_{lu}$ and $a_{lv}$ is nonzero.
  \item
    For any node $i\in\N_G(l)\setminus\{u,v\}$, $a_{li}=0$.
  \end{enumerate}
  Then the inequality $\vct{a}^\trans\vct{x}\le a_0$ is not facet
  inducing for $\CutP(G)$ unless it is a triangle inequality on
  $l,u,v$.
\end{lemma}

\begin{proof}
  The proof is by contradiction.
  Suppose the inequality $\vct{a}^\trans\vct{x}\le a_0$ is facet
  inducing for $\CutP(G)$ but it is not a triangle inequality on
  $l,u,v$.

  First we consider the case where $a_{lu}=-\lambda\le0$ and
  $a_{lv}=-\mu\le0$.
  Without loss of generality, we assume that $\lambda\le\mu$.
  Then the inequality
  \begin{equation}
    \lambda x_{uv}-\lambda x_{lu}-\mu x_{lv}
    =\lambda(x_{uv}-x_{lu}-x_{lv})-(\mu-\lambda)x_{lv}\le0
    \label{eq:degree2-triangle}
  \end{equation}
  is valid for $\CutP(G)$.
  By assumption~(ii), $\lambda$ and $\mu$ are not both zero, and the
  left hand side of the inequality~(\ref{eq:degree2-triangle}) is not
  identically zero.

  The inequality $\vct{a}^\trans\vct{x}\le a_0$ is the sum of
  (\ref{eq:degree2-triangle}) and an inequality
  \begin{equation}
    \sum_{ij\in E\setminus\{lu,lv,uv\}}a_{ij}x_{ij}
    +(a_{uv}-\lambda)x_{uv}\le a_0.
    \label{eq:degree2-rest}
  \end{equation}
  By assumption~(iii), the node $l$ is not used in the
  inequality~(\ref{eq:degree2-rest}).
  The inequality~(\ref{eq:degree2-rest}) comes from the inequality
  $\vct{a}^\trans\vct{x}\le a_0$ by collapsing the node $l$ to the
  node $v$, and is therefore valid for $\CutP(G)$.
  Therefore, the inequality $\vct{a}^\trans\vct{x}\le a_0$ is a sum of
  two valid inequalities.
  By our assumption that the inequality $\vct{a}^\trans\vct{x}\le a_0$
  is facet inducing for $\CutP(G)$, the
  inequality~(\ref{eq:degree2-rest}) is identically zero (especially
  $a_0=0$) and the inequality~(\ref{eq:degree2-triangle}) is facet
  inducing for $\CutP(G)$.
  The inequality~(\ref{eq:degree2-triangle}) is facet inducing only if
  $\lambda=\mu$, and if this holds, then the
  inequality~(\ref{eq:degree2-triangle}) is a triangle inequality on
  $l,u,v$.
  This means that $\vct{a}^\trans\vct{x}\le a_0$ is the triangle
  inequality.
  This contradicts our assumption.

  Now we consider the cases where at least one of $a_{lu}$ or $a_{lv}$
  is positive.  Switching the inequality on an appropriate subset of
  $\{u,v\}$, we can make both $a_{lu}$ and $a_{lv}$ nonpositive.
  This reduces general cases to the case where $a_{lu}\le0$ and
  $a_{lv}\le0$ hold.
\end{proof}

\begin{proof}[Proof of Proposition~\ref{prop:not-triangle}]
  Suppose the contrary: the inequality $\vct{a}^\trans\vct{x}\le a_0$
  is completely supported by the edge set
  $\{u_il,v_il\mid l\in\N_G(u_i)\cap\N_G(v_i)\}$.
  Since $G(\vct{a})$ has more than three nodes, there exists a node
  $l\in U_i\setminus\{u_i,v_i\}$ such that at least one of $a_{u_il}$
  and $a_{v_il}$ is nonzero.
  By Lemma~\ref{lemma:degree2} with $u=u_i$ and $v=v_i$, the
  inequality $\vct{a}^\trans\vct{x}\le a_0$ is not facet inducing for
  $\CutP(G)$, a contradiction.
\end{proof}

\section{Triangular elimination from $\K_n$}
  \label{sect:k-partite}

Triangular elimination from the complete graph to another graph is
useful because much is known about facets of the cut polytope of the
complete graph.

\subsection{Facets and triangular elimination from $\K_n$}
  \label{subsect:facet-kmn}

Theorem~\ref{theorem:trielim} provides a sufficient condition for an
inequality obtained by triangular elimination to be facet inducing.
We prove another sufficient condition when $G$ is the complete graph.

\begin{theorem} \label{theorem:multiple-times}
  Let $G=(V,E)$ be the complete graph on $n$ nodes with $n\ge5$.
  Let $V=V_1\cup\dots\cup V_m$ be a partition of $V$ to $m$ disjoint
  sets of nodes.
  We denote by $E_l=\{u_{l1}v_{l1},\dots,u_{lt_l}v_{lt_l}\}$ the set
  of edges in the clique on $V_l$, where
  $t_l=\abs{E_l}=\binom{\abs{V_l}}{2}$.
  Let $F=E_1\cup\dots\cup E_m$.
  Let $G'=(V',E')$ be a graph with
  $n+\sum_{1\le l\le m}t_l$ nodes.
  $n$ nodes in $G'$ are labelled by $V$, and we group the other nodes
  into $m$ sets $W_1,\dots,W_m$ with $\abs{W_l}=t_l$.
  We denote the nodes in $W_l$ by $w_{l1},\dots,w_{lt_l}$.
  If the following conditions apply, then $G'$ is a triangular
  elimination of $G$ with respect to $F$ associating node $w_{li}$
  with edge $u_{li}v_{li}$, and the triangular
  elimination of any non-triangle facet inducing inequality for
  $\CutP(G)$ is facet inducing.
  \begin{enumerate}[(i)]
  \item
    The subgraph of $G'$ induced by $V$ is the complete $m$-partite
    graph $\K_{\abs{V_1},\dots,\abs{V_m}}$ whose nodes are partitioned
    as $V_1,\dots,V_m$.
  \item
    For $l=1,\dots,m$, $W_l$ is an independent set in $G'$.
  \item
    For $l=1,\dots,m$ and $i=1,\dots,t_l$,
    $u_{li}w_{li},v_{li}w_{li}\in E'$.
  \end{enumerate}
\end{theorem}

\begin{proof}
  By conditions~(i) and (iii), it is straightforward to check that
  $G'$ is a triangular elimination of $G$ with respect to $F$.

  Let $\vct{a}^\trans\vct{x}\le a_0$ be a facet inducing inequality
  of $\CutP(G)$ which is not the triangle inequality, and
  $(\vct{a}')^\trans\vct{x}\le a'_0$ be a triangular elimination of
  $\vct{a}^\trans\vct{x}\le a_0$.
  Similarly to the proof of Theorem~\ref{theorem:trielim}, without
  loss of generality, we assume $a_0=a'_0=0$ and that no triangle
  forms used in the process of the triangular elimination has a
  nonzero constant term.

  The idea is to apply Theorem~\ref{theorem:trielim} $m$ times to
  convert $\vct{a}^\trans\vct{x}\le0$ of $\CutP(G)$ to
  $(\vct{a}')^\trans\vct{x}\le0$ of $\CutP(G^{(m)})$ where $G'$ is
  a subgraph of $G^{(m)}$, and then project the resulting facet to a
  facet of $\CutP(G')$ by using Lemma~\ref{lemma:projection}.

  First we define intermediate graphs $G^{(l)}=(V^{(l)},E^{(l)})$ and
  inequalities $(\vct{a}^{(l)})^\trans\vct{x}\le0$ for
  $l=0,1,\dots,m$.
  Let $G^{(0)}=G$ and $\vct{a}^{(0)}=\vct{a}$.
  For $l=1,\dots,m$, $G^{(l)}=(V^{(l)},E^{(l)})$ is defined by
  $V^{(l)}=V^{(l-1)}\cup W_l$ and
  $E^{(l)}=(E^{(l-1)}\setminus E_l)\cup
          \{vw\mid v\in V^{(l-1)},w\in W_l\}$.
  Then $G^{(l)}$ is a triangular elimination of $G^{(l-1)}$ with
  respect to $E_l$ where node $w_{li}\in W_l$ of $G^{(l)}$ is
  associated with edge $u_{li}v_{li}\in E_l$ of $G^{(l-1)}$.

  Let $(\vct{a}^{(l)})^\trans\vct{x}\le0$ be a triangular elimination
  of $(\vct{a}^{(l-1)})^\trans\vct{x}\le0$.

  $(\vct{a}^{(0)})^\trans\vct{x}\le0$ is facet inducing for
  $\CutP(G^{(0)})$, and the support graph of
  $(\vct{a}^{(0)})^\trans\vct{x}\le0$ has more than three nodes.
  Since triangular elimination never decreases the number of nodes of
  the support graph of an inequality, the support graph of
  $(\vct{a}^{(l)})^\trans\vct{x}\le0$ has more than three nodes for
  $l=1,\dots,m$.
  By applying Proposition~\ref{prop:not-triangle} and
  Theorem~\ref{theorem:trielim} $m$ times,
  $(\vct{a}^{(m)})^\trans\vct{x}\le0$ is facet inducing for
  $\CutP(G^{(m)})$.

  $(\vct{a}^{(m)})^\trans\vct{x}\le0$ is a triangular elimination of
  $\vct{a}^\trans\vct{x}\le0$ since
  \[
    (\vct{a}^{(m)})^\trans\vct{x}-\vct{a}^\trans\vct{x}
    =\sum_{1\le l\le m}
      ((\vct{a}^{(l)})^\trans\vct{x}-(\vct{a}^{(l-1)})^\trans\vct{x})
  \]
  is the sum of the triangular forms used in $m$ applications of
  triangular elimination.
  This combined with Proposition~\ref{prop:switching} implies that
  $(\vct{a}^{(m)})^\trans\vct{x}\le0$ is switching equivalent to
  $(\vct{a}')^\trans\vct{x}\le0$.

  By conditions~(i) and (ii), $G'$ is a subgraph of $G^{(m)}$.
  The support graph of $(\vct{a}^{(m)})^\trans\vct{x}\le0$ is a
  subgraph of a graph $G''=(V',E'')$ obtained from
  $\K_{\abs{V_1},\dots,\abs{V_m}}$ by adding nodes in
  $W_1\cup\dots\cup W_m$ and edges $u_{li}w_{li},v_{li}w_{li}$ for
  $l=1,\dots,m$ and $i=1,\dots,t_l$.
  By conditions~(i) and (iii), this support graph is a subgraph of
  $G'$.
  By Lemma~\ref{lemma:projection}, the dimension of the face of
  $\CutP(G')$ defined by $(\vct{a}^{(m)})^\trans\vct{x}\le0$ is
  $\abs{E'}-1$, which implies that
  $(\vct{a}^{(m)})^\trans\vct{x}\le0$ is facet inducing for
  $\CutP(G')$.
  Since $(\vct{a}^{(m)})^\trans\vct{x}\le0$ is switching equivalent to
  $(\vct{a}')^\trans\vct{x}\le0$, $(\vct{a}')^\trans\vct{x}\le0$ is
  also facet inducing for $\CutP(G')$.
\end{proof}

\begin{remark}
  If $a_{u_{li}v_{li}}=0$ for some edge $u_{li}v_{li}\in F$, then the
  associated node $w_{li}$ is not
  used in the triangular elimination
  $(\vct{a}')^\trans\vct{x}\le a'_0$, and the triangular elimination
  becomes facet inducing for $\CutP(G'-w_{li})$, where
  $G'-w_{li}$ denotes a graph obtained by removing node
  $w_{li}$ and edges incident to it from $G'$.
\end{remark}

\begin{corollary} \label{cor:k-partite}
  Let $G=(V,E)$, $V_l$, $E_l$, $F$, $W_l$ and $V'$ as stated in
  Theorem~\ref{theorem:multiple-times}.
  We partition $V'$ into $k$ ($m\le k\le 2m$) disjoint sets
  $V'_1,\dots,V'_k$, and let $G'=(V',E')$ be the complete $k$-partite
  graph with vertices partitioned into the sets $V'_1,\dots,V'_k$.
  If the following conditions are satisfied, then $G'$ is a triangular
  elimination of $G$, and the triangular elimination of any
  non-triangle facet inducing inequality for $\CutP(G)$ is facet
  inducing.
  \begin{enumerate}[(i)]
  \item
    For $l=1,\dots,m$, $V_l$ and $W_l$ are completely contained in some
    $V'_i$ and $V'_j$, respectively, and $i\ne j$.
  \item
    For $1\le l<l'\le m$, $V_l$ and $V_{l'}$ are contained in different
    sets $V'_i$ and $V'_j$ ($i\ne j$).
  \end{enumerate}
\end{corollary}

Theorem~2.1 of \cite{AviImaItoSas-JPA05} is the special case of
Corollary~\ref{cor:k-partite} with $m=k=3$, $\abs{V_3}=1$ (which
implies $W_3=\varnothing$), and $G'$ is the complete tripartite graph
with nodes partitioned into three sets $V_1\cup W_2$, $V_2\cup W_1$
and $V_3$, except that
Theorem~2.1 of \cite{AviImaItoSas-JPA05} also deals with the
triangular elimination of the triangle inequality.

\subsection{Triangular elimination from $\K_n$ to $\K_{r,s}$ and
  equivalence of inequalities}

Here we focus on the case $m=k=2$ in Corollary~\ref{cor:k-partite},
and we consider how Proposition~\ref{prop:switching} extends to
include permutation equivalence of inequalities.
Before that, we restate Corollary~\ref{cor:k-partite} in this case.

\begin{corollary} \label{cor:bipartite}
  Let $G=(V,E)$ be the complete graph on $n=p+q\ge5$ nodes
  $\A_1,\dots,\A_p,\allowbreak\B_1,\dots,\B_q$, and let $G'=(V',E')$ be the
  complete bipartite graph $\K_{r,s}$ with
  $r=p+\binom{q}{2},s=q+\binom{p}{2}$ where the nodes are partitioned
  into
  $\{\A_i\mid 1\le i\le p\}\cup\{\A_{jj'}\mid 1\le j<j'\le q\}$ and
  $\{\B_j\mid 1\le j\le q\}\cup\{\B_{ii'}\mid 1\le i<i'\le p\}$.
  Then $G'$ is a triangular elimination of $G$ with respect to
  $F=\{\A_i\A_{i'}\mid 1\le i<i'\le p\}\cup
     \{\B_j\B_{j'}\mid 1\le j<j'\le q\}$,
  and the triangular elimination of any non-triangle facet inducing
  inequality for $\CutP(G)$ is facet inducing.
\end{corollary}

Even if two facet inducing inequalities $\vct{a}^\trans\vct{x}\le a_0$
and $\vct{b}^\trans\vct{x}\le b_0$ of $\CutP(G)$ are equivalent up to
permutation and switching, their triangular eliminations to
$\CutP(G')$ are generally not, since different edges in $G$ may be
treated in different ways in the course of triangular elimination.
However, if we consider triangular elimination from $\CutP(\K_n)$ to
$\CutP(\K_{r,s})$ as described in Corollary~\ref{cor:bipartite},
then we know exactly when the triangular eliminations of
$\vct{a}^\trans\vct{x}\le a_0$ and $\vct{b}^\trans\vct{x}\le b_0$ are
equivalent up to permutation and switching.

\begin{theorem} \label{theorem:equivalent-kmn}
  Let $n=p+q\ge5$, $r=p+\binom{q}{2}$ and $s=q+\binom{p}{2}$, and
  label the nodes of $\K_n$ and $\K_{r,s}$ as described in
  Corollary~\ref{cor:bipartite}.
  For two non-triangle facet inducing inequalities
  $\vct{a}^\trans\vct{x}\le a_0$ and $\vct{b}^\trans\vct{x}\le b_0$
  of $\CutP(\K_n)$ and their respective triangular eliminations
  $(\vct{a}')^\trans\vct{x}\le a'_0$ and
  $(\vct{b}')^\trans\vct{x}\le b'_0$ to $\CutP(\K_{r,s})$, the
  following two conditions are equivalent.
  \begin{enumerate}[(a)]
  \item
    The two inequalities $\vct{a}^\trans\vct{x}\le a_0$ and
    $\vct{b}^\trans\vct{x}\le b_0$ can be transformed to each other by
    applying some combination of
 the switching operation, the permutation operation within
    $\{\A_1,\dots,\A_p\}$ or within $\{\B_1,\dots,\B_q\}$, and if
    $p=q$, the permutation operation swapping $A_i$ and $B_i$ for all $i$.
  \item
    The two inequalities $(\vct{a}')^\trans\vct{x}\le a'_0$ and
    $(\vct{b}')^\trans\vct{x}\le b'_0$ are permutation-switching
    equivalent.
  \end{enumerate}
\end{theorem}

\begin{remark}
  Condition~(a) implies condition~(b) even if
  $\vct{a}^\trans\vct{x}\le a_0$ and $\vct{b}^\trans\vct{x}\le b_0$
  are the triangle inequality, but the converse does not hold.
  Here is a counterexample: let $p=2$ and $q=3$.
  Let $\vct{a}^\trans\vct{x}\le a_0$ and
  $\vct{b}^\trans\vct{x}\le b_0$ be the triangle inequalities
  $-x_{\A_1\A_2}-x_{\A_1\B_1}+x_{\A_2\B_1}\le0$ and
  $-x_{\A_1\B_1}+x_{\A_1\B_2}-x_{\B_1\B_2}\le0$, respectively, and
  $(\vct{a}')^\trans\vct{x}\le a'_0$ and
  $(\vct{b}')^\trans\vct{x}\le b'_0$ be inequalities
  $-x_{\A_1\B_1}+x_{\A_2\B_1}-x_{\A_1\B_{12}}-x_{\A_2\B_{12}}\le0$
  and
  $-x_{\A_1\B_1}+x_{\A_1\B_2}-x_{\A_{12}\B_1}-x_{\A_{12}\B_2}\le0$.
  Note that the condition~(a) does not hold since $p\ne q$, whereas
  $(\vct{a}')^\trans\vct{x}\le a'_0$ and
  $(\vct{b}')^\trans\vct{x}\le b'_0$ are permutation-switching
  equivalent and the condition~(b) is satisfied.
\end{remark}

Now we give a proof of Theorem~\ref{theorem:equivalent-kmn}.

\begin{proof}
  First we prove $\text{(a)}\implies\text{(b)}$.
  If $\vct{a}^\trans\vct{x}\le a_0$ and
  $\vct{b}^\trans\vct{x}\le b_0$ are switchings of each other, then
  their triangular eliminations are also switching of each other by
  Proposition~\ref{prop:switching}.

  If $\vct{a}^\trans\vct{x}\le a_0$ is transformed to
  $\vct{b}^\trans\vct{x}\le b_0$ by swapping $\A_i$ and $\A_{i'}$,
  then the triangular elimination of $\vct{a}^\trans\vct{x}\le a_0$
  is transformed to the triangular elimination of
  $\vct{b}^\trans\vct{x}\le b_0$ by swapping $\A_i$ and $\A_{i'}$
  and swapping $\B_{ii''}$ and $\B_{i'i''}$ for all $i''\ne i,i'$,
  if we also apply this permutation to $\Tri_r$ for $1\le r\le t$.

  If $p=q$ and $\vct{a}^\trans\vct{x}\le a_0$ is transformed to
  $\vct{b}^\trans\vct{x}\le b_0$ by swapping $\A_i$ and $\B_i$ for
  $1\le i\le p$ at the same time, then the triangular elimination of
  $\vct{a}^\trans\vct{x}\le a_0$ is transformed to the triangular
  elimination of $\vct{b}^\trans\vct{x}\le b_0$ by swapping $\A_i$
  and $\B_i$ for $1\le i\le p$ and swapping $\A_{ii'}$ and
  $\B_{ii'}$ for $1\le i<i'\le p$ at the same time.

  Next we prove $\text{(b)}\implies\text{(a)}$.
  Let $(\vct{a}')^\trans\vct{x}\le a'_0$ and
  $(\vct{b}')^\trans\vct{x}\le b'_0$ be the triangular eliminations
  of $\vct{a}^\trans\vct{x}\le a_0$ and
  $\vct{b}^\trans\vct{x}\le b_0$ from $\CutP(\K_n)$ to
  $\CutP(\K_{r,s})$, respectively.
  We require that the triangle inequality $\Tri(u_i,v_i,w_i)$ was
  not used in triangular elimination to produce these two
  inequalities.
  Then $a'_0=a_0$ and $b'_0=b_0$.
  In addition, this requirement guarantees
  $a_{\A_i\A_{i'}}=\max\{a'_{\A_i\B_{ii'}},a'_{\A_{i'}\B_{ii'}}\}$
  and
  $a_{\B_j\B_{j'}}=\max\{a'_{\A_{jj'}\B_j},a'_{\A_{jj'}\B_{j'}}\}$,
  and similar equations for the vectors $\vct{b}$ and $\vct{b}'$.
  By Proposition~\ref{prop:switching}, we only need to consider the
  case where $a_0=b_0=0$ and $(\vct{a}')^\trans\vct{x}\le a_0$ and
  $(\vct{b}')^\trans\vct{x}\le b_0$ are equivalent up to
  permutation.

  The key to proving the assertion $\text{(b)}\implies\text{(a)}$
  is that from Lemma~\ref{lemma:degree2}, we can
  distinguish the nodes $\A_i$ from the nodes $\A_{jj'}$ by examining the
  inequality $(\vct{a}')^\trans\vct{x}\le0$.

  We prove the assertion~(a) holds by case analysis on the
  permutation used to transform $(\vct{a}')^\trans\vct{x}\le0$
  to $(\vct{b}')^\trans\vct{x}\le0$.
  The automorphism group of $\K_{r,s}$ is generated by permutations
  within $\{\A_1,\dots,\A_p,\allowbreak\A_{12},\dots,\A_{q-1,q}\}$,
  permutations within
  $\{\B_1,\dots,\B_q,\allowbreak\B_{12},\dots,\B_{p-1,p}\}$, and if
  $r=s$, the permutation $\tau_0$ which swaps $\A_i$ and $\B_i$ for
  $1\le i\le p$ and swaps $\A_{ii'}$ and $\B_{ii'}$ for
  $1\le i<i'\le p$ at the same time.
  Since $p+q\ge5$, $r=s$ if and only if $p=q$.

  If $p=q$ and $(\vct{a}')^\trans\vct{x}\le0$ is transformed to
  $(\vct{b}')^\trans\vct{x}\le0$ by the permutation $\tau_0$,
  then $\vct{a}^\trans\vct{x}\le0$ is transformed to
  $\vct{b}^\trans\vct{x}\le0$ by swapping $\A_i$ and $\B_i$ for
  $1\le i\le p$ at the same time.
  Therefore, from now on, we can assume that
  $(\vct{a}')^\trans\vct{x}\le0$ is transformed to
  $(\vct{b}')^\trans\vct{x}\le0$ by a permutation $\tau$ which
  permutes nodes within
  $\{\A_1,\dots,\A_p,\A_{12},\dots,\A_{q-1,q}\}$ and nodes within
  $\{\B_1,\dots,\B_q,\B_{12},\dots,\B_{p-1,p}\}$.

  Recall that the support graph of a vector $\vct{a}$ is a subgraph
  $G(\vct{a})=(V(\vct{a}),E(\vct{a}))$ of $G$ with all the edges $e$
  in $G$ with $a_e\ne0$ as its edges and all the endpoints of edges
  in $E(\vct{a})$ as its nodes.
  From Lemma~\ref{lemma:degree2}, all the nodes in $G(\vct{a})$ have
  a degree more than two.
  Since triangular elimination does not change the degree of
  existing nodes in the support graph, the nodes $\A_i$ and $\B_j$,
  if present in $G(\vct{a}')$, have degree more than two in
  $G(\vct{a}')$.
  On the other hand, from the definition of triangular elimination,
  the nodes $\A_{jj'}$ and $\B_{ii'}$, if present in $G(\vct{a}')$,
  have degree equal to two.
  Therefore, we can partition the nodes of $\K_{r,s}$ into
  three groups: $V_1$ consists of those which do not appear in
  $G(\vct{a}')$, $V_2$ consists of those with degree equal to two,
  and $V_3$ consists of those with degree more than two.
  The nodes $\A_i$ belong to $V_1$ or $V_3$, and the nodes
  $\A_{jj'}$ belong to $V_1$ or $V_2$.
  The same argument applies to $G(\vct{b}')$, and we partition the
  nodes of $\K_{r,s}$ into $W_1$, $W_2$ and $W_3$ in a parallel way.
  The permutation $\tau$ maps $V_1$ to $W_1$, $V_2$ to $W_2$ and
  $V_3$ to $W_3$, respectively.
  We define a permutation $\sigma$ on
  $\{\A_1,\dots,\A_p,\B_1,\dots,\B_q\}$ as follows.
  If $\A_i\in V_3$, then let $\sigma(\A_i)=\tau(\A_i)$.
  If $\B_j\in V_3$, then let $\sigma(\B_j)=\tau(\B_j)$.
  The rest of $\sigma$ is defined so that $\sigma$ maps the nodes in
  $V_1$ of the form $\A_i$ to the nodes in $W_1$ of the form $\A_i$,
  and the nodes in $V_1$ of the form $\B_j$ to the nodes in $W_1$ of
  the form $\B_j$.

  We show that $\sigma$ maps $\vct{a}^\trans\vct{x}\le0$ to
  $\vct{b}^\trans\vct{x}\le0$.
  All we have to prove is that for all edges $uv$ in $\K_n$, we have
  $a_{uv}=b_{\sigma(u)\sigma(v)}$.
  If $u$ belongs to $V_1$, then $\sigma(u)\in W_1$, and we have
  $a_{uv}=b_{\sigma(u)\sigma(v)}=0$.
  Since the same applies for the case $v\in V_1$, we only need to
  consider the case where both $u$ and $v$ belongs to $V_3$.
  In this case, $\sigma(u)=\tau(u)$ and $\sigma(v)=\tau(v)$.
  If $u=\A_i$ and $v=\B_j$, then
  $a_{\A_i\B_j}=a'_{\A_i\B_j}=b'_{\tau(\A_i)\tau(\B_j)}
  =b_{\tau(\A_i)\tau(\B_j)}=b_{\sigma(\A_i)\sigma(\B_j)}$.
  If $u=\A_i$ and $v=\A_{i'}$, then
  $a_{\A_i\A_{i'}}=\max\{a'_{\A_i\B_{ii'}},a'_{\A_{i'}\B_{ii'}}\}
   =\max\{b'_{\tau(\A_i)\tau(\B_{ii'})},b'_{\tau(\A_{i'})\tau(\B_{ii'})}\}
   =b_{\tau(\A_i)\tau(\A_{i'})}=b_{\sigma(\A_i)\sigma(\A_{i'})}$.
  The same applies to the case where $u=\B_j$ and $v=\B_{j'}$.
  Therefore, $\vct{a}^\trans\vct{x}\le0$ is transformed to
  $\vct{b}^\trans\vct{x}\le0$ by the permutation $\sigma$.
\end{proof}

\section{Concluding remarks}
  \label{sect:conclusion}

Theorems~\ref{theorem:trielim} and \ref{theorem:multiple-times} are
sufficient conditions for a triangular elimination of a facet inducing
inequality to be facet inducing.
An open problem is: what are necessary and sufficient conditions on
graphs $G$ and $G'$ for a triangular elimination of a non-triangle
facet inducing
inequality to be facet inducing?
Extending Theorem~\ref{theorem:equivalent-kmn} to general graphs is
another open problem.

\section*{Acknowledgments}

The third author is supported by the Grant-in-Aid for JSPS Fellows.
We would like to thank an anonymous referee for pointing out several
shortcomings in the original manuscript.

\bibliography{bell-full,bell}

\end{document}